\documentclass[12pt,reqno]{article}
mm \textheight=8.9in

\usepackage{amssymb}
\usepackage{amsmath}
\usepackage{amsthm}
\usepackage{pb-diagram}
\usepackage{color}

\newcommand{\br}[3]{{$#1$}$\lower4pt\hbox{$\tp\atop\raise4pt \hbox{$\scriptscriptstyle{#2}$}$} ${$#3$}}
\newcommand{\tw}[3]{{$#1$}${\,\scriptscriptstyle {#2}}\atop\raise9pt\hbox{$\scriptstyle\tp$} ${$#3$}}
\newcommand{\ttps}[2]{{#1}\raise5pt\hbox{$\lower12pt\hbox{$\scriptstyle\tp$}\atop \lower0pt\hbox{$\tilde\;$}$}\raise4.5pt\hbox{${\scriptstyle{#2}}$}}
\newcommand{\st}[1]{\mbox{${\,\scriptscriptstyle {#1}}\atop\raise5.5pt\hbox{$*$}$}}

\newcommand{\rd}[1]{\mbox{${\,\scriptscriptstyle {#1}}\atop\raise5.5pt\hbox{$\bullet$}$}}
\newcommand{\rt}[1]{\otimes_\chi}
\newcommand{\lt}[1]{\mbox{${\,\scriptscriptstyle {#1}}\atop\raise5.5pt\hbox{$\ltimes$}$}}
\newcommand{\btr}{\raise1.2pt\hbox{$\scriptstyle\blacktriangleright$}\hspace{2pt}}
\newcommand{\btl}{\raise1.2pt\hbox{$\scriptstyle\blacktriangleleft$}\hspace{2pt}}
\newcommand{\rl}[1]{\stackrel{#1}{\rightharpoondown}}

\newcommand{\lcr}{\raise1.0pt \hbox{${\scriptstyle\rightharpoonup}$}}
\newcommand{\rcr}{\raise1.0pt \hbox{${\scriptstyle\leftharpoonup}$}}

\newcommand{\ttp}{{\lower12pt\hbox{$\tp$}\atop \hbox{$\tilde\;$}}}

\newcommand{\id}{\mathrm{id}}

\newcommand{\Ha}{{H}}
\newcommand{\A}{{A}}
\newcommand{\B}{{B}}
\newcommand{\Lc}{{L}}
\newcommand{\J}{J}

\newcommand{\D}{\mathfrak{D}}

\newcommand{\Mc}{\mathcal{M}}

\newcommand{\Z}{\mathbb{Z}}

\newcommand{\tp}{\otimes}
\newcommand{\vt}{\vartheta}

\newcommand{\U}{U}

\newcommand{\F}{F}
\newcommand{\ve}{\varepsilon}
\newcommand{\gm}{\gamma}

\newcommand{\la}{\lambda}
\newcommand{\tr}{\triangleright}
\newcommand{\tl}{\triangleleft}

\newcommand{\End}{\mathrm{End}}

\newcommand{\Span}{\mathrm{Span}}
\newcommand{\Aut}{\mathrm{Aut}}

\newcommand{\Hom}{\mathrm{Hom}}

\newcommand{\Ann}{\mathrm{Ann}}

\newcommand{\Ga}{\Gamma}

\newcommand{\h}{\mathfrak{h}}

\newcommand{\nn}{\nonumber}

\newcommand{\al}{\alpha}

\newcommand{\be}{\begin{eqnarray}}
\newcommand{\ee}{\end{eqnarray}}

\newtheorem{thm}{Theorem}[section]
\newtheorem{propn}[thm]{Proposition}
\newtheorem{lemma}[thm]{Lemma}
\newtheorem{corollary}[thm]{Corollary}

\theoremstyle{definition}
\newtheorem{remark}[thm]{Remark}
\newtheorem{definition}[thm]{Definition}
\newtheorem{example}[thm]{Example}

\newcount\prg

\newcommand{\parag}{\advance\prg by1 {\noindent\bf\thesection.\the\prg\hspace{6pt}}}

\begin{document}
\title{On dynamical smash product}
\author{A. I. Mudrov}
\date{}
\maketitle
\begin{center}
{
 St.-Petersburg Department of Steklov Mathematical Institute,
\\
Fontanka 27, 191023 St.-Petersburg, Russia,\\
and\\
Department of Mathematics, University of York, YO10 5DD, UK }\\
 \texttt{e-mail: mudrov@pdmi.ras.ru}
\end{center}
\begin{abstract}
In the theory of dynamical Yang-Baxter equation, with any Hopf algebra $\Ha$ and a certain
$\Ha$-module and $\Ha$-comodule algebra $\Lc$ (base algebra) one associates a monoidal category.
Given an algebra $\A$ in that category, one can construct an associative
algebra $\A\rtimes \Lc$, which is a generalization of the ordinary smash product
when $\A$ is an ordinary $\Ha$-algebra. We study this "dynamical smash product"
and its modules  induced from one-dimensional representation of
the subalgebra $\Lc$. In particular, we construct an analog
of the Galois map $\A\tp_{\A^\Ha} \A\to \A\tp \Ha^*$.
\end{abstract}
\section{Introduction}
This work is motivated by a recently discovered connection between equivariant deformation quantization and
the dynamical Yang-Baxter equation (YBE), see \cite{DM1,AL,KMST,EEM}. One goal of this paper
is to give a deeper algebraic insight to constructions arising along those lines,
in particular, to non-associative algebras participating in quantization, \cite{DM1}.
To a certain extent, this work is also about non-commutative "principal bundles".
Such a point of view was already suggested in \cite{DM1}, and here we develop that idea one step further.

Conventionally, the function algebra of a non-commutative  principal fiber bundle is understood as  Hopf-Galois
extension of the subalgebra of invariants.
 What is special about our point of view as compared to the standard one is
that we do not assume the principal bundle to be an algebra. The necessity of
such a generalization is dictated by the fact
that for a great deal of interesting examples the principal bundle cannot be quantized in the class of
associative algebras.
Rather, we consider functions on a principal bundle as sections of the associated
vector bundle whose fiber is the function algebra on the structure group (viewed as the co-regular module). In this sense
it can be quantized for important homogeneous spaces of simple groups. Thus we
view non-commutative principal bundles as modules over non-commutative quotient spaces
subject to a Galois-like condition, which we introduce and study in this paper.

The term "dynamical" takes its origin in the mathematical physics literature,
where the corresponding generalization of the YBE appears in connection with integrable
conformal field theories. Contrary to the purely algebraic ordinary YBE,
its dynamical analog is differential in the classical  version and difference-like in the quantum.
The theory of dynamical YBE provides a natural environment for deformation quantization
of $G$-spaces alongside with vector bundles over them. In this geometrical interpretation the set of dynamical variables (more precisely, its invariant
classical points) appears to be the moduli space
of quantized principal bundles. This moduli space admits an action of the "Picard group" of the associated linear bundles.

The ground object of our theory is the so called dynamical base,
which is the pair $(\Ha,\Lc)$ of a Hopf algebra $\Ha$ and a certain $\Ha$-algebra $\Lc$.
Algebraically, $\Lc$ is a braided-commutative algebra in the quasitensor category of modules over the double of $\Ha$.
Our present interest is centered on a dynamical algebra over $(\Ha,\Lc)$,
which is an algebra in a certain monoidal category rather than in the usual sense. In the
simplest non-trivial case $\Lc=\Ha$,  the dynamical algebras
can be described as follows. Consider an associative algebra $\B$ containing $\Ha$ as a subalgebra.
We endow $\B$ with the structure of  left $\Ha$-module under the adjoint action. Suppose
there exists an $\Ha$-submodule  $\A\subset \B$  freely generating $\B$ over $\Ha$.
Then $\A$ is a dynamical algebra and $\B$ is a generalization of the smash product
of $\A$ and $\Ha$. A similar description is valid for general $\Lc$, and
we reserve  for $\B$ the notation $\A\rtimes \Lc$.

For every dynamical base $(\Ha,\Lc)$  we define and study a homomorphism onto a base $(\Ha',\Lc')$
with commutative $\Ha'$ and $\Lc'$, which reduction can be accompanied with
a reduction of a dynamical algebra $\A$; the reduced dynamical algebra
is a certain submodule in $\A$. We study the reduced case in detail and classify
all reduced dynamical bases. This can serve to a rough classification of dynamical bases
and dynamical algebras.

With any $\Ha$-invariant character $\chi$ of the algebra $\Lc$ we relate
an associative algebra $\A_\chi^\Lc$ and a pair of right and left $\A_\chi^\Lc$-modules
$\A_\chi$ and ${}_{\tilde\chi}\A$ on the same underlying vector space $\A$. Here $\tilde\chi$ is another
character canonically constructed out of $\chi$.
These modules play the role of non-commutative left and right "principal bundles".
They are constructed as $\A\rtimes \Lc$-modules via right and left induction from the one dimensional representations $\chi$
and $\tilde \chi$ of the subalgebra $\Lc\subset \A\rtimes \Lc$;
then $\A_\chi^\Lc$ is realized as the algebra of $\A\rtimes \Lc$-endomorphisms
of  $\A_\chi$ and ${}_{\tilde\chi} \A$. This is the algebraic meaning
of $\A_\chi^\Lc$, which appeared in \cite{DM1} as the deformed function algebra on
a coadjoint orbit (conjugacy class) of a simple complex group. In the special case
of \cite{DM1} the vector space $\A_\chi^\Lc$ coincides with $\Ha$-invariants in $\A$.

For every quadruple $(\Ha,\Lc,\A,\chi)$ we
construct a $\B$-bimodule $\A\rt{\chi}\Ha^*$ on the tensor product $\A\tp \Ha^*$.
Here $\Ha^*$ is the dual Hopf algebra to $\Ha$.
We introduce a $\B$-bimodule map $\A_\chi\tp _{\A_\chi^\Ha}{}_{\tilde\chi}\A\to \A\rt{\chi}\Ha^*$,
where $\A_\chi^\Ha\subseteq \A^\Lc_\chi$ stands for the subalgebra of invariants.
This is an analog of the Galois map for dynamical algebras and it turns into
the ordinary Galois map when $\A$ is an ordinary $\Ha$-algebra. We call the pair $(\A,\A^\Ha_\chi)$
Galois extension if this map is an  isomorphism.
Restriction of the dynamical Galois map to $\A^\Lc_\chi\tp _{\A_\chi^\Ha}\A^\Lc_\chi$
coincides with the usual Galois map for the extension $\A^\Lc_\chi/{\A_\chi^\Ha}$.

Finally, for every invariant character $\chi$
in the left regular representation we construct
a Morita context between  $\A\rtimes \Ha$ and $\A^\Ha_\chi$.
We construct it also between $\A\rtimes \Lc$ and  $\A^\Lc_\chi$,
out of a $\chi$-generator $e\in \Lc$ in the left regular representation of $\Lc$.

Throughout the paper $\Ha$ is assumed to be a finite dimensional Hopf algebra over a field $\Bbbk$.
Some of  our results can be generalized to a certain class of infinite dimensional Hopf algebras that are close
to quantum groups, but we make this simplifying assumption by technical reasons in order to avoid
numerous stipulations.

\vspace{0.2in}
\noindent
{\bf \large Acknowledgements}
This research is supported by the EPSRC grant C511166 and partly
supported by the RFBR grant 06-01-00451. The author is grateful
to the Max-Planck Institute f\"{u}r Mathematik in Bonn for hospitality.

\section{Preliminaries (dynamical base)}
Fix a finite dimensional Hopf algebra $\Ha$ over a field $\Bbbk$ with the
comultiplication $\Delta$, the counit $\ve$ and the antipode $\gm$.
The antipode is invertible in a finite dimensional Hopf algebra, and we denote its inverse by $\bar\gm$.
We shall use the same notation for the structure operations in the dual Hopf algebra $\Ha^*$; this
will not cause any confusion, as elements of $\Ha$ and $\Ha^*$ will be always explicitly stated.

We denote by $\D(\Ha)$ the double Hopf algebra of $\Ha$ and we choose its realization
in the form $\Ha\bowtie \Ha_{op}^*$, where $op$ designates the opposite multiplication.
In other words, the double is constructed on the tensor product of $\Ha\tp\Ha_{op}^*$,
where the factors  are embedded as sub-bialgebras. Given a basis $\{h_i\}\subset \Ha$ and its dual
$\{\eta^i\}\subset \Ha^*$, the element $\Theta:=\sum_i \eta^i\tp h_i\in  \D(\Ha)$ is
a quasitriangular structure on $\D(\Ha)$; it is independent on the choice of basis.
The Yang-Baxter equation on $\Theta$ is equivalent to the commutation relations between
elements of $\Ha$ and of $\Ha^*_{op}$ and thus determines the algebra structure of $\D(\Ha)$.

The symbol $\tp$ stands for the tensor product over the ground field $\Bbbk$;
the tensor product over other rings will be indicated explicitly.
We adopt the following Sweedler-like convention denoting symbolically the coproducts and coactions:
$$
\Delta(h)=h^{(1)}\tp h^{(2)},
\quad
\Delta(\eta)=\eta^{\langle 1\rangle }\tp \eta^{\langle 2 \rangle },
$$
for $h\in \Ha$ and $\eta\in \Ha^*$. Thus we use
different indication for the Sweedler $\Ha$- and $\Ha^*$-components.
For a left  $\Ha$-comodule $V$  and a right $\Ha^*$-comodule $W$ we use the symbolic presentation
$$
\delta(v)=v^{(1)}\tp v^{[\infty]},
\quad
\delta(w)=w^{[0]}\tp w^{\langle 1\rangle},
$$
where $v\in V$ and $w\in W$.
We also suppress summation when labelling tensor factors of tensor objects
and denote them simply as $\Theta=\Theta_1\tp \Theta_2$ for $\Theta\in \D(\Ha)^{\tp 2}$ {\em etc}.

We denote by $\rcr$ and $\lcr$ the right and, respectively, left co-regular actions of $\Ha$ on $\Ha^*$:
Explicitly they are given by
\be
h\lcr \eta= \eta^{\langle 1\rangle}\langle\eta^{\langle 2\rangle},h \rangle,
\quad \eta\rcr h= \langle\eta^{\langle 1\rangle},h \rangle\eta^{\langle 2\rangle},
\label{coreg}
\ee
for $h\in \Ha$ and $\eta\in \Ha^*$. Here $\langle ., .\rangle$ designates the
canonical pairing between $\Ha^*$ and $\Ha$.

Recall that a Yetter-Drinfeld  (YD) module $V$ over $\Ha$ is a left $\Ha$-module and left $\Ha$-comodule obeying
$$
\delta(h\tr v)=h^{(1)}v^{(1)}\gm(h^{(3)})\tp h^{(2)}\tr v^{[\infty]}
$$
for all $v\in V$ and $h\in \Ha$;
the symbol $\tr$ stands for the action of $\Ha$ on $V$.
Such modules form a quasitensor category, which for finite dimensional
$\Ha$ is isomorphic to the category of $\D(\Ha)$-modules.

Recall that an $\Ha$-module algebra $\B$ is an algebra in the monoidal category of $\Ha$-modules.
In other words, the multiplication $\B\tp \B\to \B$ is supposed to be equivariant.
Similarly, an $\Ha$-comodule algebra is an algebra in the category of $\Ha$-comodules.
That means that the coaction is an algebra homomorphism.
An algebra in the YD category is simultaneously a YD module, $\Ha$-module algebra
and $\Ha$-comodule algebra.  In the finite dimensional case under consideration, a YD
algebra is a module algebra over $\D(\Ha)$.

Every YD module $V$ is equipped with an $\Ha$-equivariant permutation $\tau_{V,W}\colon V\tp W\to W\tp V$ with
every $\Ha$-module $W$.
The permutation acts by the assignment $\tau_{V,W}(v\tp w)= v^{(1)}\tr w\tp v^{[\infty]}$ for all $v\in V$ and $w\in W$.
A YD-algebra is called braided-commutative if its multiplication is stable  under this permutation.

The following definition is the key ingredient of the theory of dynamical Yang-Baxter equation.
\begin{definition}
A {\em dynamical base} (or simply base) is a pair $(\Ha,\Lc)$ of a Hopf algebra $\Ha$ and a
braided commutative  YD algebra $\Lc$.
\end{definition}
We call $\Lc$ an $\Ha$-base algebra or simply a base algebra when $\Ha$ is clear from the context.
We always assume $\Lc$ to be unital.
Here we write down, for the reader's convenience, the defining axiom of base algebra $\Lc$:
\begin{enumerate}
\item $\Lc$ is a left $\Ha$-module algebra under the action $\tr\colon \Ha\tp \Lc\to \Lc$,
\item $\Lc$ is a left $\Ha$-comodule algebra under the coaction
$
\delta\colon \Lc\to \Ha\tp \Lc$.

\item the action and coaction satisfy the YD condition
$
\delta(h\tr \la)=h^{(1)}\la^{(1)}\gm(h^{(3)})\tp h^{(2)}\tr \la^{[\infty]}
$
for all $h\in \Ha$, $\la\in \Lc$.
\item the braided commutativity condition $\la\mu=(\la^{(1)}\tr\mu) \la^{[\infty]}$ for all $\la,\mu\in \Lc$ holds true for all $\la,\mu\in \Lc$.
\end{enumerate}
The above definition works for arbitrary Hopf algebras.
for finite dimensional $\Ha$
a base algebra $\Lc$ is precisely a module algebra
over the double $\D(\Ha)$ satisfying the braided commutativity condition.
In terms of the universal R-matrix  $\Theta$, the $\Ha$-coaction can be written as
\be
\label{coactions}
\la\mapsto \Theta_2\tp \Theta_1\tr\la =\la^{(1)}\tp \la^{[\infty]},
\ee
for all $\la\in \Lc$. The  braided commutativity is then reads
$\la\mu=(\Theta_2\tr \mu)(\Theta_1\tr \la)$, where the $\Ha$-action $\tr$ is extended to
the action of the double $\D(\Ha)$.

We also consider the base algebra $\Lc$ as a right comodule algebra over $\Ha^*$ corresponding to
the left $\Ha$-action by duality. The same applies to every $\Ha$-module (algebra).

A trivial example of base algebra is given by the ground field $\Bbbk$, considered as the trivial $\Ha$-module and $\Ha$-comodule.
Below we give two more examples of base algebras of which the first is valid not only for finite dimensional but for arbitrary
 Hopf algebras.
\begin{example}
\label{L=H}
Consider $\Ha$ as a left $\Ha$-module algebra with respect to the adjoint action and a left $\Ha$-comodule algebra
under the comultiplication. Then the conditions 3-4 are satisfied, and $(\Ha,\Ha)$ is a dynamical base.
\end{example}
\begin{example}
\label{L=H*}
Let $\{h_i\}\subset \Ha$ be a basis and $\{\eta^i\}\subset \Ha^*$ be its dual basis.
Consider $\Ha^*$ as a left  $\Ha$-module algebra with respect to the
action (\ref{coreg}) and a left $\Ha$-comodule algebra with respect
to the coaction $a\mapsto \sum_{i,j}\gm(h_i)h_j\tp \eta^ia \eta^j$; this definition does not depend on the choice of basis.
Then conditions 3-4 are satisfied, and $(\Ha,\Ha^*)$ is a dynamical base.
\end{example}
\begin{definition}
A homomorphism $(\Ha_1,\Lc_1)\to (\Ha_2,\Lc_2)$ of dynamical bases
is a pair $\varphi\colon \Ha_1\to \Ha_2$, $\varpi\colon \Lc_1\to \Lc_2$ of
mappings, where
$\varphi$ is a Hopf algebra homomorphism and $\varpi$ an algebra homomorphism
respecting the actions and coactions:
\be
\label{base_homo}
\varpi(h\tr \la)=\varphi(h)\tr\varpi(\la),
\quad (\varphi\tp \varpi)\bigl(\delta_1(\la)\bigr)=\delta_2 \bigl(\varpi(\la)\bigr),
\ee
for all $h\in \Ha_1$ and $\la\in \Lc_1$.
\end{definition}
\begin{example}
\label{al-auto}
Consider $\Ha$ as an $\Ha$-base algebra, as an Example \ref{L=H}.
For any group-like element $\al\in \Ha^*$ the map $\varpi_\al\colon h\mapsto \al\lcr h=h^{(1)}\langle\al,h^{(2)}\rangle$
is an algebra automorphism. It is easy to check that the pair $(\id, \varpi_\al)$ an automorphism of dynamical
base $(\Ha,\Ha)$.
\end{example}

Suppose that  $(\Ha_1,\Lc_1)$ is a dynamical base, $\Lc_2$ is a YD algebra over $\Ha_2$,
and the pair of homomorphisms $(\varphi,\varpi)$ obeys (\ref{base_homo}). If $\varpi$ is onto,
then $\Lc_2$ is braided commutative, and hence $(\Ha_2,\Lc_2)$ is a dynamical base.
Then the pair $(\varphi,\varpi)$ is a homomorphism of dynamical bases.

\section{Invariant characters of base algebra}
\label{SecICh}
In the present section we study invariant characters of base algebras, by which we understand
unital $\Ha$-equivariant homomorphisms of $\Lc$ to the ground field regarded as the trivial $\Ha$-module. We are especially interested
in the situation when the corresponding one-dimensional representation is realized as a left ideal in $\Lc$, e.g.
when $\Lc$ is finite dimensional semisimple.

Consider a unital homomorphism  $\chi\colon \Lc\to \Bbbk$ of algebras.
We call it invariant character if $\chi(h\tr \la)=\ve(h)\chi(\la)$ for all $h\in \Ha$ and $\la\in \Lc$. Note that
for $\Lc=\Ha$ from Example \ref{L=H} any character is invariant. On the contrary, the $\Ha$-base algebra $\Ha^*$
from Example \ref{L=H*} has no invariant
characters, unless $\dim\Ha=1$.
Denote by $\hat\Lc$ the set of invariant characters of $\Lc$. The set $\hat\Ha\subset \Ha^*$  is a group,
with the identity $\ve$ and the inverse $\bar\al=\al\circ \gm=\al\circ \bar\gm$; here $\gm$ is the antipode in $\Ha$.
The group $\hat\Ha$ naturally acts on $\Lc$ on the right by
$\la\cdot \al=(\al\tp \id)\circ \delta (\la)$. This action induces a left action on $\hat \Lc$ by
$(\al\cdot \chi)(\la)=(\al\tp \chi)\circ \delta (\la)$.

Every invariant character defines  a homomorphism
$\iota_\chi\colon\Lc\to \Ha$ through the mapping
$\la\mapsto \la^{(1)}\chi(\la^{[\infty]})$ for all $\la\in \Lc$.
It amounts to a homomorphism of $\Ha$-bases $(\Ha,\Lc)\to (\Ha,\Ha)$ identical on $\Ha$.

 For any invariant character $\chi\in \hat\Lc$ we call $e\in \Lc$ a left $\chi$-generator
if  $\Lc e=\chi(\Lc) e$.
Let $\Lc^\chi$ denote the vector space of left $\chi$-generators.
In the special case $\Lc=\Ha$ and $\chi=\ve$, it coincides with the space $\int_l$ of left integrals.
By construction, $\Lc^\chi$ is a two-sided ideal in $\Lc$.
\begin{lemma}
\label{gen-inv}
The vector space $\Lc^\chi$ is $\Ha$-invariant.
\end{lemma}
\begin{proof}
Indeed, for all $h\in \Ha$, $\la\in \Lc$, and $e\in \Lc^\chi$:
$$
\la (h\tr e)=\bigl(h^{(2)}\bar\gm(h^{(1)})\tr\la\bigr) (h^{(3)}\tr e)=h^{(2)}\tr\bigl((\bar\gm(h^{(1)})\tr\la)e\bigr)
=\chi(\la)\ve\bigl(\bar\gm(h^{(1)}\bigr)h^{(2)}\tr e.
$$
The last expression is equal to $\chi(\la)h\tr e$.
\end{proof}

We call a character $\chi\in \hat\Lc$ projective if the corresponding one dimensional
$\Lc$-module
is projective. In other words, if $\Lc^\chi\not\subset \ker\chi$ or, equivalently,
if $\Lc^\chi$ contains an idempotent. Below we show that $\dim \Lc^\chi=1$
if the
 character $\chi$ is  projective.

Next we introduce an operation on invariant characters through  a Drinfeld element in the double of $\Ha$,
which is defined for any finite dimensional quasitriangular Hopf algebra, \cite{Dr2}.
Set
$$
\vt=\Theta_1\gm^{-2}(\Theta_2),
\quad
\bar\vt=\Theta_1\gm(\Theta_2).
$$
It is known from \cite{Dr2} that $\bar \vt=\vt^{-1}$. Conjugation with $\vt$ implements the squared antipode in $\D(\Ha)$:
$\vt h \vt^{-1}=\gm^2(h)$.
The element $\vt$ satisfies the equation
\be
\label{gauge}
\Theta_{21}\Theta_{12}\Delta(\vt)=(\vt\tp \vt).
\ee
If follows from  (\ref{gauge}) that the action of $\vt$ on $\Lc$
implements an algebra automorphism $\la\mapsto \vt\tr \la$. That is a consequence
of braided-commutativity of $\Lc$, so that $\Theta$-s in the left hand side of (\ref{gauge})
vanish on $\Lc$. The automorphism $\vt$ commutes with
the action of $\hat \Ha\subset \D(\Ha)$, as the square antipode is identical on $\hat\Ha$.

For every $\chi\in \hat \Lc$ define the adjoint character $\tilde \chi$ by
the assignment $\tilde\chi\colon \la \mapsto \chi(\bar\vt\tr\la)$. It is invariant by virtue of
$\vt h=\gm^2(h)\vt$ for all $h\in \Ha$.
Since the action of $\vt$ commutes with the action of $\hat \Ha$, we have
 $\al\cdot\tilde\chi=\widetilde{\al\cdot\chi}$ for all $\al\in \hat \Ha$.
We call a character self adjoint if $\tilde \chi=\chi$.
Next we characterize the
mapping $\chi\to \tilde \chi$ in terms of a certain action of $\Lc$ on $\Ha^*$.

Define a left $\Lc$-action on $\Ha^*$ by setting
$
\la\rl{\chi}\eta:= \eta\rcr(\gm\circ\iota_\chi)(\la)
$
for all $\la\in \Lc$ and $\eta\in \Ha^*$.

\begin{lemma}
\label{scalar}
For all $\la\in \Lc$ the element $\la^{[0]}\rl{\chi}\la^{\langle 1\rangle}\in \Ha^*$
is the  scalar multiple $\tilde \chi(\la)1$.
\end{lemma}
\begin{proof}
Consider $\Lc$ as an algebra over the double $\D(\Ha)$.
Evaluating $\la^{[0]}\rl{\chi}\la^{\langle 1\rangle}$ at an arbitrary element $h\in \Ha$ we find
$$
\langle \la^{[0]}\rl{\chi}\la^{\langle 1\rangle},h\rangle=\chi(\Theta_{1} \Theta_{2'}\tr\la) \langle \Theta_{1'},\gm(\Theta_{2}) h\rangle=
\chi(\Theta_{1}\gm(\Theta_{2}) h \tr\la)=\chi(\bar\vt h \tr\la)=\chi(\bar \vt\tr\la)\ve(h).
$$
Here we use the dashed subscripts to distinguish between different copies of $\Theta$.
The last equality is obtained by pulling $\bar\vt$ to the right
and taking into account  $\bar \vt h \vt =\bar\gm^2(h)$.
\end{proof}

Similarly to left $\chi$-generators one can  define $\chi$-generators under the right regular representation
of $\Lc$: those are elements $e\in \Lc$ satisfying  $e\Lc=e\chi(\Lc)$.
\begin{propn}
\label{left-right}
For all $\chi\in \hat \Lc$, left $\chi$-generators are right $\tilde \chi$-generators and {\em vice versa}.
\end{propn}
\begin{proof}
Let $e\in \Lc^\chi$ be a left $\chi$-generator. Then, for all $\la\in \Lc$,
$$
e \la = (e^{(1)}\tr \la) e^{[\infty]}=\langle \la^{\langle 1\rangle},e^{(1)} \rangle \la^{[0]} e^{[\infty]}
=\langle \la^{[0]}\rl{\chi}\la^{\langle 1\rangle} ,e^{(1)}\rangle e^{[\infty]}
=\tilde \chi(\la) \ve(e^{(1)}) e^{[\infty]}=\tilde \chi(\la)e.
$$
In the middle equality we have used the identity $e^{(1)}\tp  \mu e^{[\infty]}=(\gm\circ\iota_\chi)(\mu) e^{(1)}\tp  e^{[\infty]}$
for $\mu=\la^{[0]}\in \Lc$. The next equality to the right is due to Lemma \ref{scalar}.
This calculation proves that $e$ is a right $\tilde\chi$-generator. To check that
every right $\tilde\chi$-generator is a left $\chi$-generator,
notice that the opposite algebra $\Lc_{op}$ is a base algebra
over $\Ha_{op}$, cf. \cite{DM2}.
The subscript $op$ designates the opposite multiplication. The action of $\Ha_{op}$ on $\Lc_{op}$
is given by the assignment $h\tp \la\mapsto \bar\gm(h)\tr\la$, and the coaction stays the same
as the $\Ha$-coaction on $\Lc$.
This readily implies that every right $\tilde \chi$-generator is a left $\chi$-generator.
\end{proof}

\section{Base reduction}
In the present section we  study a  reduction procedure for dynamical base.
Namely, with any  base $(\Ha,\Lc)$ such that $\hat \Lc\not= \emptyset$ we associate a dynamical base $(\Ha',\Lc')$
with commutative $\Ha'$ and $\Lc'$ and a base homomorphism $(\Ha,\Lc)\stackrel{\phi,\varpi}{\longrightarrow}(\Ha',\Lc')$,
which is a surjective on the components. Then we classify  reduced dynamical bases.
\begin{propn}
\label{1-dim}
Suppose an invariant character $\chi\in \hat \Lc$ is projective. Then a) the
$\Ha$-module $\Lc^\chi$ is trivial and contained
in the center of $\Lc$,
b) $\dim \Lc^\chi=1$, c) $\chi$ is self-adjoint.
\end{propn}
\begin{proof}
Suppose there exists a generator $e\in \Lc^\chi$ such that $\chi(e)\not=0$.
We can assume that $e$ is an idempotent and hence $\chi(e)=\tilde\chi(e)=1$, due to Proposition \ref{left-right}.
If  $e'\in \Lc^\chi$ is yet another $\chi$-generator, then  the equalities
$$
\chi(e+e')e=(e+e')e=(e+e')\tilde\chi (e)
$$
imply $\chi(e')e=e'$. This proves b).

First of all remark that every $\Ha$-invariant in $\Lc$ are central.
Now c) follows from a) and Proposition \ref{left-right}.
To prove a) recall from Lemma \ref{gen-inv} that $\Lc^\chi=\Bbbk e$ is an $\Ha$-module. This module is trivial, as the character $\chi$ is $\Ha$-invariant.
This proves a).
\end{proof}
\begin{remark}
Remark that Propositions \ref{left-right} and \ref{1-dim} generalize  the well known fact that
the vector space of left integrals $\int_l$ in a finite dimensional Hopf algebra is one dimensional, \cite{Sw}.
On setting $\Lc=\Ha$ an integral $t\in \int_l$ becomes  a left $\ve$-generator
a right $\tilde \ve$-generator. However $\dim \int_l=1$ always
for $\Lc=\Ha$, contrary to general base algebra case.
For example, consider the algebra $\Lc=\Span\{1,e_i\}_{i=1}^k$ with the nil multiplication on $\Span\{e_i\}_{i=1}^k$.
Endowed with the structure of trivial $\D(\Ha)$-module, $\Lc$  becomes an $\Ha$-base algebra.
The character $\chi(1)=1$, $\chi(e_i)=0$ is invariant and $\dim \Lc^\chi=k$,
because $\Lc^\chi$ in this example coincides with $\ker \chi$.
\end{remark}

Consider the subspace $\hat\Ha^\bot$ in $\Ha$ annihilated by all characters of $\Ha$,
that is, the intersection $\hat\Ha^\bot=\cap_{\al\in \hat\Ha} \ker \al$.
It is a bi-ideal, and the quotient $\Ha/\hat\Ha^\bot$ is a Hopf algebra. The latter is commutative,
and can be viewed as the function algebra $\Bbbk^{\hat \Ha}$ on the group $\hat \Ha$. Put $\Ha':=\Ha/\hat\Ha^\bot$ and denote by $\varpi$
the corresponding projection $\Ha\to\Ha'$.

From now one we suppose that the set of invariant characters of $\Lc$ is not empty.
Let $\hat\Lc^\bot$ be the intersection  $\cap_{\chi\in \hat\Lc}\ker \chi$, which is an ideal
in $\Lc$. Set $\Lc':=\Lc/\hat\Lc^\bot$
and denote by $\varphi$ the projection $\Lc\to\Lc'$.
As the characters are unital homomorphisms and $\hat \Lc\not= \emptyset$ by assumption,
one has $\Lc'\not=\{0\}$. Clearly the algebra $\Lc'$ is commutative.
\begin{propn}
The pair $(\Ha',\Lc')$ is a dynamical base,
and
$(\Ha,\Lc)\stackrel{\varphi,\varpi}{\longrightarrow}(\Ha',\Lc')$ is a homomorphism
of bases.
\end{propn}
\begin{proof}
The ideal $\hat \Lc^\bot$ is obviously $\Ha$-invariant. Moreover,
$\hat \Ha^\bot \Lc\subset \Lc^\bot$. Indeed, for all $h\in \hat \Ha^\bot$ and $\la\in \Lc$
one has $\chi(h\tr \la)=\ve(h)\chi(\la)=0$. Hence $\Lc'=\Lc/\hat \Lc^\bot$
is a natural module algebra over $\Ha'$. By construction,
$\varpi(h\tr \la)=\varphi(h)\varpi(\la)$, that is, the projection $\varpi$ is $\varphi$-equivariant.
 Under the coaction $(\varphi\tp \id)\circ \delta$, the ideal $\hat \Lc^\bot$ is also a coideal.
Thus, $\Lc'$ is a YD-module algebra over $\Ha'$ and hence a base algebra.
This proves the assertion.
\end{proof}
\begin{definition}
A dynamical base $(\Ha,\Lc)$ is called reduced if $\cap_{\chi\in \hat \Lc}\ker\chi=\{0\}$.
The reduced base $(\Ha',\Lc')$ constructed above is called reduction of $(\Ha,\Lc)$.
\end{definition}

If $\Lc'$ is finite dimensional, then it is semisimple. Otherwise, we can select in $\Lc'$ a finite dimensional subalgebra,
which is also a base algebra over $\Ha'$.
\begin{propn}
Suppose that the base $(\Ha,\Lc)$ is reduced   and $\oplus_{\chi\in \hat \Lc}\Lc^\chi\not=\{0\}$.
Then $\oplus_{\chi\in \hat \Lc}\Lc^\chi$ is a base algebra over $\Ha$.
\end{propn}
\begin{proof}
It is sufficient to show that $\oplus_{\chi}\Lc^\chi$ is a sub-comodule over $\Ha'$.
Suppose $\Lc^\chi\not =0$ and let $e_\chi\in \Lc^\chi$ be the idempotent generating $\Lc^\chi$, cf. Proposition \ref{1-dim}.
For all $\al\in\hat \Ha$ the element $\al(e_\chi^{(1)})e_\chi^{[\infty]}$ is the idempotent
$e_{\bar\al\cdot \chi}$ generating $\Lc^{\bar\al\cdot \chi}$. This readily implies the statement.
\end{proof}
Note that $\oplus_{\chi}\Lc^\chi$
is not only a subalgebra of $\Lc$ but also a quotient algebra.
Indeed, $\oplus_{\chi}\Lc^\chi$ is the image under the homomorphism $\la\mapsto \sum_{\chi\in \hat \Lc}e_\chi \la$.

Any finite dimensional reduced base algebra may be regarded as a function algebra on the
set $\hat\Lc$. The action of the group  $\hat \Ha$ on $\hat\Lc$ stratifies $\hat \Lc$ to a disjoint union of
orbits. Restriction of $\Lc$ to every orbit is again a base algebra over $\Ha$; the whole
$\Lc$ splits into a direct sum of these base algebras.

Below we prove a statement, which gives rise to classification of semisimple reduced dynamical bases.
We consider a slightly more general situation: a) the group
$\hat \Ha^*$ acts transitively on a set $X$, b) $\Lc$ is the function algebra
$\Bbbk^{X}$, c) the characters of $\Bbbk^{X}$, which are
naturally identified with elements of $X$, are not necessarily invariant.
\begin{thm}
Let $G$ be a finite group, $K\subset G$ its subgroup, and $X=G/K$ the set of left  cosets.
The base algebra structures on $\Bbbk^{X}$ are parameterized by central idempotents $\pi\in \Bbbk K$.
The base $(G,K,\pi)$ is reduced if and only if $\pi$ is the identity, otherwise the
set $\widehat{\Bbbk^{X}}$  of invariant characters is empty.
Two bases $(G,K_i,\pi_i)$ are isomorphic if and only if there exits an automorphism
$\varphi\in \Aut(G)$ such that $K_2=\varphi(K_1)$ and $\pi_2=\varphi(\pi_1)$.
\end{thm}
\begin{proof}
The left action of $G$ on $X$ corresponds to the right action  of $\Bbbk G$ on $\Bbbk^{X}$,
which gives rise to a left $\Bbbk^{G}$-coaction.
Let $e_\psi\in \Lc$ denote the idempotent corresponding to a point $\psi\in G/K$.
The commutative algebra $\Bbbk^{G}$ acts on $\Lc$ by $h\tr e_\psi=\langle\pi(\psi),h\rangle e_\psi$
for some $\pi(\psi)\in \Bbbk G$. The condition $e_\psi e_\chi=\delta_{\psi\chi}e_\chi$, with $\delta_{\psi\chi}$ being the
Kronekker symbol,
implies that $\Lc$ is a $\Bbbk^{G}$-module algebra if and only if $\pi(\psi)$ is an idempotent in $\Bbbk G$.
This idempotent  is the group identity if and only if the character $\psi$ is invariant.

Further, the YD condition is equivalent to $G$-equivariance of the map $\pi\colon X\to G$,  that is, to the equality
$\pi(\al\cdot\psi)=\al\pi(\psi)\al^{-1}$ for all $\al\in G$ and $\psi\in X$. Hence $\pi(\psi)$ is determined by its value
$\pi(\psi_0)$ where $\psi_0\in X$ is the origin. The idempotent $\pi(\psi_0)$ commutes with $K$, the isotropy subgroup of $\psi_0$.

The inclusion $\pi(\psi_0)\in \Bbbk K$ is equivalent to  $\D(\Bbbk^{G})$-commutativity.
Indeed, the equation $e_\psi e_\phi=(e_\psi^{(1)}\tr e_\phi) e_\psi^{[\infty]}$ is turns into
$\phi(e_\psi) e_\phi=\langle \pi(\phi),e_\psi^{(1)}\rangle \phi(e_\psi^{[\infty]}) e_\phi $
for all $\phi,\psi\in X$. This is equivalent to
$\phi(e_\psi) =\bigl(\pi(\phi)\cdot\phi \bigr)(e_\psi)$ for all  $\phi,\psi\in X$.
That is, to
$\phi =\pi(\phi)\cdot\phi $ for all  $\phi\in X$. This is the case
if and only if $\pi(\psi_0)\in K$.

The part concerning the isomorphism is an easy exercise left to the reader.
\end{proof}
\section{The bialgebra $\Ha^\chi$}
In the present section we study certain properties of $\Ha$-modules in the presence of
the base algebra $\Lc$. We shall heavily rely on these technical results in the sequel.
As another application, we associate a bialgebra $\Ha^\chi$ with every invariant character $\chi\in \hat \Lc$.
This bialgebra is constructed as a quotient of $\Ha$ and also is used in what follows.

Let $V$ be a  left $\Ha$-module.
For every integer $k$ and ordered pair $(\psi,\chi)\in \hat \Lc\times \hat \Lc$ we define $V[k,\psi,\chi]\subset V$
to be the vector subspace
of $\psi$-generators in $V$ under the representation
$\Lc\stackrel{\iota_\chi}{\longrightarrow}\Ha\stackrel{\gm^{k}}{\longrightarrow}\Ha\to \End(V)$,
that is,
$v\in V[\psi,\chi]\Leftrightarrow (\gm^{k}\circ\iota_\chi)(\la)v=\psi(\la)v,
$
for all $\la\in \Lc$.
Note that this representation is left if $k$ is even and right if $k$ is odd.
It is easy to check that $V[k,\psi,\chi]$ is an $\Ha$-submodule in $V$.
For every $\al\in \hat\Ha$ and all $\chi\in \Lc$ the $\Ha$-submodule of weight $\al$ is
obviously contained in $V[k,\chi,\al^{(-1)^k}\cdot\chi]$.
\begin{lemma}
\label{l-r}
For all $k\in \Z$ one has
$$V[k,\psi,\chi]=V[k\pm 1,\chi,\psi],\quad V[k,\psi,\chi]=V[k,\tilde\psi,\tilde \chi].$$
\end{lemma}
\begin{proof}
For all $v\in V[k,\psi,\chi]$ and all $\la\in \Lc$ we find
$$
\begin{array}{ccc}
\gm^{k\pm 1}(\la^{(1)})\chi(\la^{[\infty]}) v =\gm^{k\pm 1}(\la^{(1)})\gm^{k}(\la^{(2)})\psi(\la^{[\infty]}) v
=\psi(\la) v,
\\[5pt]
\chi(\la^{[\infty]})\gm^{k\mp 1}(\la^{(1)}) v =\psi(\la^{[\infty]})\gm^{k}(\la^{(2)})\gm^{k\mp 1}(\la^{(1)}) v
=\psi(\la) v,
\end{array}
$$
where we take the upper sign  if $k$ is even and the lower sign if $k$ is odd.
The left equality in the second line is obtained using the fact that $V[k,\psi,\chi]$ is $\Ha$-invariant.
The above calculation shows that  $V[k,\psi,\chi]\subset V[k\pm 1,\chi,\psi]$ for all $k\in \Z$.
This implies the first equality of the proposition.

The right equality follows from the left one, because for all $\la$ the operator
$$
\tilde\psi(\la^{(\infty)}) \gm^{k}(\la^{(1)}) =\psi(\bar\vt\tr  \la^{(\infty)}) \gm^{k}(\la^{(1)}) =
\psi\bigl((\bar\vt\tr  \la)^{(\infty)}\bigr)  \gm^{k+2}\bigl((\bar\vt\tr  \la)^{(1)}\bigr)
$$
is the scalar multiple $\chi(\bar\vt\tr  \la)=\tilde \chi(\la)$ on $V[k,\psi,\chi]=V[k+2,\psi,\chi]$.
The right equality is obtained using the presentation (\ref{coactions}) for
the coaction, and the standard fact $(\gm\tp \gm)(\Theta)=\Theta$.
The above calculation proves the inclusion $V[k,\tilde\psi,\tilde\chi]=V[k,\psi,\chi]$. The inverse inclusion is
proved similarly.
\end{proof}
The above proposition implies  that the module $V[k,\psi,\chi]$ actually depends on $k \mod 2$, and
$V[k,\chi,\chi]$ is independent on $k$ at all.
In what follows, we reserve the notation $V[\psi,\chi]:=V[0,\psi,\chi]$ and $V[\chi]:=V[\chi,\chi]$.
\begin{remark}
\label{left-right_antipode}
One can equally consider right $\Ha$-modules and define the corresponding right and left
representations of $\Lc$. The statement analogous to Lemma \ref{l-r} will be true as
well. That can be reduced to the considered case by passing to the dynamical base
$(\Ha_{op},\Lc_{op}$, cf. the proof of Proposition \ref{left-right}.
\end{remark}

For every pair $V,W$ of left $\Ha$ modules we find that $V[\phi,\psi]\tp W[\psi,\chi]\subset (V\tp W)[\phi,\chi]$. Therefore,
$V[\phi,\psi]\tp W[\psi,\chi]\to U[\phi,\chi]$
under every $\Ha$-equivariant mapping $V\tp W\to U$.

With every invariant character $\chi\in \hat \Lc$ we shall associate a bialgebra $\Ha^\chi$.
Let $\Mc$ be the monoidal category of left $\Ha$-modules.
Denote by $\Ann(V,\chi)\subset \Ha$ the annihilator of $V[\chi]$. Consider the
intersection $\Xi_\chi:=\cap_V \Ann(V,\chi)$ taken over all modules from $\Mc$. It is an ideal in $\Ha$,
and the inclusion $V[\chi]\tp W[\chi]\subset (V\tp W)[\chi]$ implies that it is a
bi-ideal: $\Delta(\Xi_\chi)\subset \Xi_\chi\tp \Ha+\Ha\tp \Xi_\chi$.
For $\Bbbk$ regarded as the trivial $\Ha$-module, we have $\Bbbk[\chi]=\Bbbk$, hence $\Xi_\chi$ is annihilated
by the counit. The quotient $\Ha^\chi=\Ha/\Xi_\chi$ is a bialgebra, and it is a Hopf algebra,
provided $\gm(\Xi_\chi)=\Xi_\chi$.
By construction,  for every $V$ from $\Mc$ the representation of $\Ha$ in $V[\chi]$ factors through
a representation of $\Ha^\chi$.

Recall that $\hat\Ha^\chi\subset\hat\Ha$ stands for the isotropy subgroup of the point $\chi\in \hat \Lc$.
\begin{propn}
Characters $\widehat{\Ha^\chi}$ of the bialgebra $\Ha^\chi$ constitute  a subgroup in $\hat \Ha$ isomorphic to $\hat\Ha^\chi$.
\end{propn}
\begin{proof}
As $\Ha^\chi$ is a homomorphic image of $\Ha$, there is an inclusion $\widehat{\Ha^\chi}\subset \hat\Ha$.
We have to show that $\widehat{\Ha^\chi}= \hat\Ha^\chi$.
Let $\Bbbk_\al$ denote the one dimensional representation of $\Ha$ corresponding to $\al\in \hat \Ha$.
If $\al\in \hat \Ha^\chi$, then $\Bbbk_\al[\chi]=\Bbbk_\al$. Hence $\al(\Xi_\chi)=\{0\}$ and $\al\in  \widehat{\Ha^\chi}$; this
proves the inclusion
$\hat\Ha^\chi\subset \widehat{\Ha^\chi}$.
Conversely, if $\al\in \widehat{\Ha^\chi}$, then $\al(\Xi_\chi)=\{0\}$.
For every $\la\in \Lc$ the element $\iota_\chi(\la)-\chi(\la)1\in \Ha$ belongs to $\Xi_\chi$.
Hence $0=\al\bigl(\iota_\chi(\la)\bigr)-\chi(\la)=\al\cdot \chi(\la)-\chi(\la)$ for all $\la\in \Lc$,
that is, $\al\in \hat \Ha^\chi$.
\end{proof}

Next we give a detailed consideration to a special case, when the bialgebra $\Ha^\chi$ becomes a Hopf algebra.
Namely, we assume that the character $\chi\in \hat \Lc$ is projective.
Let $e_\chi\in \hat\Lc$ be the corresponding idempotent.
The element $t_\chi:=\iota_\chi(e_\chi)\in \Ha$ is also an idempotent and it
is central in $\Ha$, because $e_\chi$ is invariant.
\begin{lemma}
For every left $\Ha$-module $V$ the submodule $V[\chi]$ has the form $V[\chi]=t_\chi V$.
\end{lemma}
\begin{proof}
Indeed, if $t_\chi v=v$ for some $v\in V$, then
$$
\iota_\chi(\la)v=\iota_\chi(\la)t_\chi v=\iota_\chi(\la e_\chi) v=\chi(\la)t_\chi v=\chi(\la) v,
$$
that is, $t_\chi V\subset  V[\chi]$.
Conversely, if $\iota_\chi(\la)v=\chi(\la)v$ for all $\la\in \Lc$, then
$$
t_\chi v=\iota_\chi(e_\chi)v=\chi(e_\chi)v=v.
$$
This proves the inverse inclusion $t_\chi V\supset  V[\chi]$.
\end{proof}
Regard $\Ha$ as a left and right $\Lc$-module via the homomorphism $\iota_\chi\colon \Lc\to \Ha$.
Then  $\Ha[\chi]$ is a subring $t_\chi\Ha t_\chi=\Ha t_\chi\subset \Ha$.
It is also isomorphic to the quotient of $\Ha$ by  the ideal $(1-t_\chi)\Ha$, which
coincides with the ideal $\Xi_\chi$ introduced above.

\begin{corollary}
The idempotent $t_\chi$ is fixed under the antipode $\gm$, and $\gm(\Ha[\chi])=\Ha[\chi]$.
\end{corollary}
\begin{proof}
For all $\la\in \Lc$ and all $h\in \Ha[\chi]$ we find
$$
\gm\bigl(\iota_\chi(\la)\bigr)h=\chi(\la)h,
\quad
h\bar\gm\bigl(\iota_\chi(\la)\bigr)=h\chi(\la).
$$
The left equality is a direct consequence of Lemma \ref{l-r}, the right one holds true
due to Remark \ref{left-right_antipode}.
Setting $h=t_\chi$ and $\la=e_\chi$ we get
$
\gm(t_\chi)t_\chi=t_\chi
$
and
$
t_\chi\bar\gm(t_\chi)=t_\chi.
$
From this we conclude that $t_\chi=\gm(t_\chi)$ and hence $\gm(\Ha[\chi])=\Ha[\chi]$.
\end{proof}

Denote by $\Ha^\chi$ the algebra $\Ha[\chi]$ and consider $t_\chi$ as the projector from
$\Ha$ to $\Ha^\chi$.
\begin{thm}
The comultiplication
\be
\label{comult_proj}
\Ha^\chi\hookrightarrow \Ha\stackrel{\Delta}{\longrightarrow}\Ha\tp \Ha \stackrel{t_\chi\tp t_\chi}{\longrightarrow}\Ha^\chi\tp \Ha^\chi
\ee
 along with the counit $\ve|_{\Ha^\chi}$ and the antipode $\gm|_{\Ha^\chi}$
makes $\Ha^\chi$ a Hopf algebra.
\end{thm}
\begin{proof}
The comultiplication is obviously an algebra homomorphism. Its coassociativity
readily follows from the equality
$$
\Delta(t_\chi)(t_\chi\tp t_\chi)=t_\chi\tp t_\chi.
$$
Let us check this equality. We have for the left-hand side:
$$
e^{(1)}_\chi t_\chi\tp \iota_\chi(e^{[\infty]}_\chi e_\chi)
=
e^{(1)}_\chi t_\chi\tp \chi(e^{[\infty]}_\chi)t_\chi
=
t_\chi t_\chi\tp t_\chi
=
t_\chi\tp t_\chi,
$$
as required.

Further, from $\ve(t_\chi)=\chi(e_\chi)=1$
we derive
$$
(\ve\tp \id)\Bigl(\Delta (h)(t_\chi\tp t_\chi)\Bigr)=1\tp ht_\chi=1\tp  h,
\quad
(\id\tp \ve)\Bigl(\Delta (h)(t_\chi\tp t_\chi)\Bigr)=ht_\chi\tp 1 =h\tp 1
\nn
$$
for all $h\in \Ha^\chi$.
Thus the restriction of $\ve$ to $\Ha^\chi$ is the counit.

To verify the antipode, notice that for all $h\in \Ha^\chi$
$$
\gm(t_\chi)\gm(h^{(1)})h^{(2)}t_\chi=\ve(h)\gm(t_\chi)t_\chi=\ve(h)t_\chi,
\quad
t_\chi h^{(1)}\gm(h^{(2)})\gm(t_\chi)=\ve(h)t_\chi\gm(t_\chi)=\ve(h)t_\chi.
\nn
$$
Hence  $\gm|_{\Ha^\chi}$ is the antipode indeed, as $t_\chi$ is the unit in  $\Ha^\chi\subset \Ha$.
This completes the proof.
\end{proof}

Note that for every character $\al\in \hat \Ha$ the bialgebras $\Ha^\chi$
and $\Ha^{\al \cdot\chi}$ are isomorphic. The isomorphism is included in the commutative diagram
$$
\begin{array}{cccc}
\Ha &\to &\Ha\\
\downarrow&&\downarrow\\
\Ha^\chi &\to &\Ha^{\al \cdot\chi}
\end{array}
$$
of bialgebra homomorphisms, where the the vertical arrows are projections and the
top arrow is given by the assignment $h\mapsto \bar\al(h^{(1)})h^{(2)}\al(h^{(3)})$, $h\in \Ha$.

\section{Dynamical algebras}

A dynamical or shifted associative algebra $\A$ over the base $(\Ha,\Lc)$ is an $\Ha$-module equipped with an equivariant map (multiplication)
$\divideontimes\colon \A\tp \A\to \A\tp \Lc$ satisfying the (shifted associativity) condition:
\be
&
\begin{diagram}
   \dgARROWLENGTH=0.5\dgARROWLENGTH
\node{\A \tp \Lc\tp \A}\arrow{e,t}{\tau_\A}\node{\A \tp \A \tp \Lc}
\arrow{e,t}{\divideontimes }\node{\A \tp \Lc \tp \Lc}\arrow{e,t}{m_\Lc}
\node{\A \tp \Lc} \arrow{s,r,!}{\parallel}
\\
\node{\A\tp\A\tp \A }\arrow{e,t}{\id\tp \divideontimes }\arrow{n,l}{\divideontimes\tp \id}
\node{\A \tp \A \tp \Lc}\arrow{e,t}{\divideontimes}\node{\A \tp \Lc \tp \Lc}\arrow{e,t}{m_\Lc}\node{\A \tp \Lc}
\end{diagram}
\label{shifted_assoc0}
\ee
Here $m_\Lc$ is the multiplication in $\Lc$ and $\tau_\A$ is the $\Ha$-invariant permutation $\Lc\tp \A\to \A\tp \Lc$ acting by $\la\tp a\mapsto \la^{(1)  }\tr a\tp \la^{[\infty]}$;
the obvious identity maps are suppressed.

The dynamical algebra is called unital if there exists an element $1_\A$ such that
$1_\A \divideontimes a=a \divideontimes 1_\A =a\tp 1_\Lc$ for all $a\in \A$.
Like an ordinary algebra,  a dynamical algebra can be made unital by extension (see the next section).

\begin{example}
Suppose $\A$ is a module algebra over a Hopf algebra.
Under its multiplication and subsequent embedding
$\A\simeq \A\tp 1_\Lc\subset \A\tp \Lc$ if becomes a dynamical algebra over an arbitrary base algebra $\Lc$.
We call such dynamical algebras  trivial.
\end{example}

\begin{example}
Suppose $\A$ is a module algebra over a Hopf algebra $\U$ containing $\Ha$
as a Hopf subalgebra.
Suppose $\F\in \U\tp \U\tp \Lc$ is a dynamical twist, \cite{DM1}.
Then the multiplication $a\divideontimes b:=(\F_1\tr a)(\F_2\tr b)\tp \F_3$
makes $\A$ a dynamical algebra.
\end{example}

Further we define a reduction procedure for dynamical algebras.
For any dynamical algebra $\A$ over the base $(\Ha,\Lc)$ it
provides a dynamical algebra $\A'$ over the reduced base $(\Ha',\Lc')$.
Consider the $\Ha$-submodule $\A'\subset \A$ annihilated by $\hat \Ha^\bot\subset \Ha$.
It is straightforward to check that the composition
$\divideontimes'\colon\A'\tp \A'\to \A\tp \Lc\to \A\tp \Lc'$ of the multiplication restricted to $\A'$ with
the projection $\varpi\colon \Lc\to \Lc'$
takes values in $\A'\tp \Lc'$. Indeed, if $h\in \hat \Ha^\bot$, then for all $\chi\in \hat\Lc$,
and $a,b\in \A'$ one has $h\tr(\id\tp \chi)(a\divideontimes b)=(\id\tp \chi)(h^{(1)}\tr a\divideontimes h^{(2)}\tr b)=0$.
Here we have used $\Ha$-invariance of the character
$\chi$, equivariance of $\divideontimes $ and the fact that $\hat \Ha^\bot$ is a bi-ideal.

\begin{propn}
The map $\A'\tp \A'\to \A'\tp \Lc'$  makes $\A'$ a dynamical algebra over the  base $(\Ha',\Lc')$.
\end{propn}
\begin{proof}
By construction, $\A'$ is a module over $\Ha'$ and the multiplication $\A'\tp \A'\to \A'\tp \Lc'$
is $\Ha'$-equivariant.
Restrict the leftmost bottom vertex in (\ref{shifted_assoc0}) to $\A'\tp \A\tp \A'$
and apply the projection $\varpi$ to the rightmost $\Lc$. This yields
the commutative diagram
\be
&
\begin{diagram}
   \dgARROWLENGTH=0.5\dgARROWLENGTH
\node{\A \tp \Lc\tp \A'}\arrow{e,t}{\tau_\A}\node{\A \tp \A' \tp \Lc}
\arrow{e,t}{\varpi}
\node{\A \tp \A' \tp \Lc'}
\arrow{e,t}{\divideontimes '}
\node{\A' \tp \Lc' \tp \Lc'}\arrow{s,t}{m_{\Lc'}}
\\
\node{\A'\tp\A'\tp \A' }\arrow{e,t}{\id\tp \divideontimes' }\arrow{n,l}{\divideontimes\tp \id}
\node{\A' \tp \A' \tp \Lc'}\arrow{e,t}{\divideontimes'}\node{\A' \tp \Lc' \tp \Lc'}\arrow{e,t}{m_{\Lc'}}\node{\A' \tp \Lc'}
\end{diagram},
\label{shifted_assoc1}
\ee
where we suppress the obvious identity morphisms.
Let $\tau'_{\A'}$ denote the canonical permutation $\Lc'\tp \A'\to \A'\tp\Lc'$.
Obviously,  $(\id\tp \varpi)\circ \tau_\A = \tau_\A '\circ  ( \varpi\tp \id)$
when restricted to $\Lc\tp \A'$. Thus we can replace all the non-dashed elements of the above diagram
by their dashed counterparts. This yields the shifted associativity diagram for
$\A'$ as a dynamical algebra over $(\Ha',\Lc')$.
\end{proof}

Remark that the reduced dynamical algebra $\A'$ contains elements
of weight $\al$ for all $\al\in \hat \Ha$. Indeed, for every such element
 $a\in \A$ and all $h\in \hat \Ha^\bot$ we find $h\tr a=\al(h) a=0$, so $a\in \A'$.

Every invariant character defines an associative algebra $I_\chi$ on the vector space $\A^\Ha$ of $\Ha$-invariants
in $\A$.
The multiplication in $I_\chi$ is set to be
$$
a\stackrel{\chi}{*}b:=(\id_\A\tp \chi)(a \divideontimes b).
$$
In applications to deformation quantization, see \cite{DM1}, $\Ha$ is a quantized universal enveloping algebra
$U(\h)$ acting on a Poisson manifold $M$. The dynamical algebra $\A$ is a "deformation" of the algebra of functions
on $M$. Then $I_\chi$ can be interpreted as a deformation of the function algebra of the corresponding
quotient space. It is not realized as a subalgebra in an associative quantization
of $M$. Rather, the function algebra on $M$ is quantized to right and left $I_\chi$-modules
$\A_\chi$ and ${}_{\tilde \chi}\A$.
Our goal is to construct a Galois map $\A_\chi\tp_{ I_\chi}{}_{\tilde \chi}\A\to \A\tp \Ha^*$ in
generalization of the standard situation when $\A_\chi={}_{\tilde \chi}\A$ is an algebra ($\A$ is a trivial dynamical algebra).
This Galois map will be also a homomorphism of bimodules over
an associative algebra with the underlying vector space $\A\tp \Lc$.
This algebra is a generalization of the ordinary smash product, and it is the subject of
our interest in the next section.
\section{Dynamical smash product}
Suppose that $\A$ is a  dynamical algebra with multiplication
$\divideontimes\colon \A\tp \A\to \A\tp \Lc$.
Define the right action of $\Lc$ on $\A\tp \Lc$ by $(a\tp \la)\mu=a\tp \la\mu$.
Introduce the associative algebra $\A\rtimes \Lc$ on the vector space $\A\tp \Lc$
with multiplication
\be
\label{dynamical_smash}
(a\tp \la)(b\tp \mu)=\Bigl(a\divideontimes(\la^{(1)}\tr b)\Bigr)(\la^{[\infty]}\mu).
\ee
The last factor is an element of $\Lc$ acting on the term within the big parentheses,
which is an element of $\A\tp \Lc$.
The algebra $\A\rtimes \Lc$ is an $\Ha$-module, and the multiplication is $\Ha$-equivariant.
\begin{definition}
The $\Ha$-algebra $\A\rtimes \Lc$ is called  dynamical smash product of $\A$ and $\Lc$.
\end{definition}
The dynamical smash product is a generalization of the usual smash product of an $\Ha$-module algebra $\A$ and
the base algebra $\Lc$,
where the dynamical multiplication $\A\tp \A\to \A\tp \Lc$ factors to the usual multiplication $\A\tp \A\to \A$
and the embedding of $\A$ in $\A\tp\Lc$.

The dynamical smash product $\B=\A\rtimes \Lc$ is also  an $\Lc$-bimodule, with
the left $\Lc$-action being defined through the right action by the formula
$\la(b\mu):=(\la^{(1)}\tr b)(\la^{[\infty]}\mu)$, for all $\la,\mu\in \Lc$ and $b\in \B$.
This can be also considered as a compatibility with the action of $\Ha$.
We call such algebras $(\Ha,\Lc)$-module algebras.
In fact, they are algebras in the monoidal category of modules over the bialgebroid $\Lc\#\Ha$, see \cite{Lu,DM2}.
Note that $\B$ is freely generated by $\A$ over $\Lc$ (with respect to either actions). These conditions characterize
 dynamical algebras.
\begin{propn}[\cite{DM1}]
Suppose $\A$ is a dynamical algebra over the base $(\Ha,\Lc)$.
Then the dynamical smash product is an $(\Ha,\Lc)$-module algebra.
Conversely, suppose $\B$ an $(\Ha,\Lc)$-module algebra.
Suppose that $\B$ is freely generated over $\Lc$ by an $\Ha$-submodule $\A\subset \B$.
Then $\A$ is a dynamical algebra.
\end{propn}

Given an $(\Ha,\Lc)$-module algebra $\B$ we can define a new $(\Ha,\Lc)$-module algebra
on the direct sum $\B\oplus \Lc$ of $\Ha$-modules.
The multiplication is given by the rule
$$
(a\oplus \la)(b\oplus \mu):=(ab+\la b+ a\mu )\oplus \la\mu.
$$
The vector $(0,1_\Lc)$ is the identity in this new algebra.
If $\B$ is a dynamical smash product $\A\rtimes \Lc$, so is the new algebra,
and the corresponding dynamical algebra $\A\oplus \Bbbk$ is unital.
From now one we assume that all the  dynamical algebras under consideration are unital.
It is then convenient to identify $\A$ with $\A\rtimes 1_\Lc\subset \B$
and  $\Lc$ with $1_A\rtimes \Lc\subset \B$.
The algebra $\B=A\rtimes \Lc$ is generated by $\A$ and $\Lc$, with the commutation relations
$h a=(h^{(1)}\tr a) h^{(2)}$.

For instance, put  $\Lc=\Ha$ and suppose that $\Ha$ is a subalgebra in $\B$.
Then $\B$ is an $\Ha$-module algebra with respect to the adjoint action.
It is also an $(\Ha,\Ha)$-module algebra. If $\B$ is freely generated over $\Ha$
by an $\Ha$-submodule $\A$, then $\A$ is a dynamical algebra.

\begin{definition}
Suppose $(\Ha,\Lc)\stackrel{\varphi,\varpi}{\longrightarrow}(\Ha',\Lc')$
is a homomorphism of dynamical bases and let $\A$ and $A'$ be dynamical algebras
over  $(\Ha,\Lc)$ and $(\Ha',\Lc')$, respectively. An equivariant mapping $\theta\colon \A\to \A'\tp \Lc'$
is called morphism of dynamical algebras if the diagram
\be
&
\begin{diagram}
   \dgARROWLENGTH=0.5\dgARROWLENGTH
\node{\A \tp \A}\arrow{s,l}{\divideontimes}\arrow{e,t}{\theta \tp\theta}\node{\A' \tp \Lc'\tp \A'\tp \Lc'}
\arrow{e,t}{\tau'_{\A'}}\node{\A' \tp \A'\tp \Lc' \tp \Lc'}\arrow{e,t}{m_{\Lc'}}
\node{\A'\tp \A'\tp \Lc'}\arrow{s,r}{\divideontimes'}
\\
\node{\A\tp\Lc}\arrow{e,t}{\theta\tp \varpi}
\node{\A' \tp \Lc' \tp \Lc'}\arrow{e,t}{m_{\Lc'}}\node{\A' \tp \Lc'}\node{\A' \tp \Lc'\tp\Lc'}\arrow{w,t}{m_{\Lc'}}
\end{diagram}
\ee
is commutative (here we have dropped the obvious identity maps).
\end{definition}

Given a homomorphism $(\Ha,\Lc)\to (\Ha',\Lc')$ of dynamical bases and their module algebras $\B$ and $\B'$ we define a homomorphism
$\B\to \B'$ as that of associative algebras which is equivariant
under the left Hopf algebra actions and the two-sided base algebra actions.
The following proposition is straightforward.
\begin{propn}
Morphisms $(\Ha,\Lc,\A)\stackrel{\varphi,\varpi,\theta}{\longrightarrow}(\Ha',\Lc',\A')$
are in natural 1-1 correspondence with equivariant homomorphisms
 $(\id\tp {m_{\Lc'}})\circ(\theta\tp \varpi)\colon\A\rtimes\Lc\to\A'\rtimes\Lc'$
of associative algebras that preserve the base algebras.
\end{propn}
\section{Linear induction from the base algebra}
\label{secLBI}
From now on we use the notation $\B$ for the dynamical smash product $\A\rtimes \Lc$ assuming it to be fixed
once and for all.
In the present section we study the modules over $\B$ that are induced
from $\Ha$-invariant characters of $\Lc$.

Denote by $\tl$ the right $\Ha$-action on $\A$ defined by
$a\tl h:= \bar \gm(h)\tr a$.
We shall use the following symbolic presentation of the dynamical multiplication:
$$
a\divideontimes b = (a\st{i}b)\rho_i=\ell_i(a\rd{i}b)=\rho_i^{[\infty]}\bigl((a\st{i}b)\tl\rho_i^{(1)}\bigr),
$$
where $a\st{i}b$ and $a\rd{i}b$ stand for the $\A$-component while
$\rho_i$ and $\ell_i$ for the $\Lc$-component. The label $i=1,2,\ldots$ is used to
distinguish between different copies of these operations rather than for summation, which
is implicitly understood as in the Sweedler notation.

For a fixed invariant character $\chi$ of the algebra $\Lc$ we regard $\A$ as the left induced $\B$-module
$\B\tp_\Lc \Bbbk$ and denote it by $\A_\chi$. Similarly, the right $\B$-module $\Bbbk\tp_\Lc \B=\A$ will be denoted by ${}_\chi\A$.
The left and right actions of $\Ha$ on $\A$ induce the left and right actions of $\Lc$
through the homomorphism $\iota_{\chi}\colon\Lc\to \Ha$.
In particular, the right $\Lc$-action on ${}_\chi\A$ is given by $a\tp \la= (\bar\gm\circ \iota_\chi)(\la)\tr
a$.
Denote by $\A_\chi^\Lc=\A_\chi[\chi]$ the space of right $\chi$-generators in $\A$ under the left $\Lc$-action
induced by the homomorphism $\iota_\chi\colon\Lc\to \Ha$.
 Due to Lemma \ref{l-r}, $\A_\chi^\Lc$ coincides
with the space  of right $\tilde\chi$-generators in the right module ${}_{\tilde \chi}\A$.

We introduce the following  $\Ha$-equivariant  bilinear operations on  $\A$:
$$
a\st{\chi}b:= a\st{1}b  \chi(\rho_1)
,\quad
a\rd{\chi} b:=\chi(\ell_1)(a\rd{1} b).
$$
The unit of the dynamical algebra $\A$  becomes the neutral element for
the operations $\st{\chi}$ and $\rd{\chi}$.
Note that $\st{\chi}$ is the restriction to $\A\subset \B$ of  the left action on $\A_\chi$;
similarly, $\rd{\chi}$ is the restriction to $\A\subset \B$ of  the right action on ${}_{\chi}\A$.

We can characterize  $\A^\Lc_\chi$ as the algebra of endomorphisms of the $\B$-modules
$\A_\chi$ and ${}_{\tilde\chi}\A$.
\begin{thm}
\label{End}
The operations $\st{\chi}$ and $\rd{\tilde\chi}$  coincide on $\A_\chi^\Lc$.
They make $\A_\chi^\Lc$
into
an associative $\Ha$-algebra  and  $\A$ into  right and left $\A_\chi^\Lc$-modules
respectively.
Moreover,
$$\End_{\B}(\A_\chi) \simeq (\A_{\chi}^\Lc)_{op}\simeq \End_{\B}({}_{\tilde\chi}\A).$$
\end{thm}
\begin{proof}
First we demonstrate that the operations in question coincide when restricted to $\A^L_\chi$.
From the presentation $\ell_1(a\rd{1} b)= \rho_1^{[\infty]}\bigl(\bar\gm(\rho_1^{(1)})\tr(a\st{1}b)\bigr)=(\Theta_1\tr \rho_1)\bigl(\bar\gm(\Theta_2)\tr(a\st{1}b)\bigr)$
 we find for all $a,b\in \A_\chi^\Lc$:
$$
a\rd{\tilde\chi}b=(a\st{1}b)\tl \iota_{\tilde\chi}(\rho_1)=
\tilde \chi \bigl(\Theta_{1}\Theta_{1'}\vt\tr \rho_1\bigr)\bigl(\bar\gm(\Theta_2)\tr a\bigr)\st{1} \bigl(\bar\gm(\Theta_{2'})\tr b)\bigr)
=(a\st{1}b) \tilde\chi\bigl(\vt\tr\rho_1\bigr)=a\st{\chi}b.
$$
The second equality follows from the $\Ha$-equivariance of the multiplication in $\B$.
In the next equality to the right we have  used Lemma \ref{l-r} allowing to drop the antipodes.

The second part of the theorem will be proved only for $\st{\chi}$, as the case of $\rd{\tilde\chi}$ is analogous.
Apply the map $(\id_\A\tp \chi)$ to the equality
$a(bc)=(ab)c$ in the algebra $\B=\A\rtimes \Lc$.
The left-hand side immediately gives
$a\st{\chi}(b \st{\chi} c)$ for all  $a,b,c\in \A$.
The right-hand side turns into $(a\st{\chi}b)\st{\chi}c$ once $c\in \A_\chi^\Lc$,
because
$a\st{1}b \tp \iota_\chi(\rho_1)\tr c=a\st{1}b \tp  \chi(\rho_1)c= (a\st{\chi}b)\tp c $
in this case.

Thus we have shown that $(a\st{\chi}b)\st{\chi}c=a\st{\chi}(b \st{\chi} c)$ whenever $c\in \A_\chi^\Lc$.
This proves that $\A_\chi^\Lc$
is an associative algebra and $\A$ is a right $\A_\chi^\Lc$-module under $\stackrel{\chi}{*}$.
Let us show that $\A_\chi^\Lc\subset \End_{\B}(\A_\chi)$. The algebra $\Lc$ acts on $\A_\chi$
by multipliers $\chi(\la)$, $\la\in \Lc$. Hence the right action of $\A_\chi^\Lc$ commutes with
the left action of $\Lc\subset \B$. The action of $\A\subset \B$ on $\A_\chi$
coincides with the operation $\stackrel{\chi}{*}$. As the operation
$\stackrel{\chi}{*}$ is associative once the third argument is from $\A^\Lc_\chi$,
the right action of $\A^\Lc_\chi$ commutes with the left action of $\A$. This gives an anti-algebra homomorphism
$\A_\chi^\Lc\to \End_{\B}(\A_\chi)$.
Note that this homomorphism is an embedding, because  $\A$ is assumed unital as a dynamical algebra and  the operation $\st{\la}$ is also unital.

The Frobenius reciprocity gives $\End_{\B}(\A_\chi)\simeq \Hom_{\Lc}(\Bbbk,\A_\chi)$.
The action of $\Lc$ on $\A_\chi$ is given by $\iota_\chi(\la)\tr a$, so
$\Hom_{\Lc}(\Bbbk,\A_\chi)=\A^\Lc_\chi$. We have shown that $\A_\chi^\Lc\simeq \End_{\B}(\A_\chi)$
as vector spaces.
It is easy to see that this is also an anti-isomorphism of algebras.
\end{proof}

Regarding $\B$ as an $\Lc$-module through the homomorphism $\iota_\chi\colon \Lc\to \Ha$,
denote by $\B^\Lc_\chi:=\B[\chi]$ the subspace of left $\chi$-generators in $\B$. By Lemma \ref{l-r},
it is at the same time  the space of   right $\tilde\chi$-generators under the right representation $\bar\gm\circ \iota_{\tilde \chi}$.
The vector space $\B^\Lc_\chi$ forms an $\Ha$-module, and equivariance of the multiplication in $\B$
implies that $\B^\Lc_\chi$ is a subalgebra in $\B$.

Introduce the left and, respectively, right $\B$-module homomorphisms
\be
\wp_\chi\colon \B\to \A_\chi,
\quad
{}_{\chi}\hspace{-1pt}\wp\colon \B\to {}_{\chi}\A,
\ee
obtained by evaluating $\chi$ on the
$\Lc$-factor in the factorizations $\A\Lc=\B=\Lc\A$.
It is also a homomorphism of $\Ha$-modules.
When restricted to $\B^\Lc_\chi$, the map $\wp_\chi$
takes values in $\A^\Lc_\chi$.
\begin{propn}
\label{BH-AH}
Restriction of $\wp_\chi$
to $\B^\Lc_\chi$ defines a surjective homomorphism $\B^\Lc_\chi\to \A^\Lc_\chi$
of $\Ha$-algebras. Through this homomorphism, $\wp_\chi\colon \B\to \A_\chi$
becomes a right $\B^\Lc_\chi$-bimodule map.
\end{propn}
\begin{proof}
As the base algebra $\Lc$ is unital, $\A^\Lc_\chi\subset \B^\Lc_\chi$ and the
$\wp_\chi$ is identical on $\A^\Lc_\chi$. Hence it is surjective.
For all $a=\sum_ia_i\la_i\in \B$ and $b=\sum_mb_m\mu_m\in \B^\Lc_\chi$ we find
$$
\wp_\chi(ab)
=\wp_\chi\Bigl(\sum_{i,m}a_i\bigl(\la_i^{(1)}\tr (b_m\mu_m)\bigr)\la_i^{[\infty]}\Bigr)
=\sum_{i,m}(a_i\st{1}b_m)\chi(\rho_1\mu_m \la_i)
=\wp_\chi(a)\st{\chi}\wp_\chi(b).
$$
In the middle equality we used the assumption  $b\in \B^\Lc_\chi$ as the latter is supposed to be invariant.
This proves two things: firstly, $\wp_\chi$ is an algebra homomorphism when restricted to $\B^\Lc_\chi$ and, secondly,
it is a homomorphism of right $\B^L_\chi$-modules $\B\to \A_\chi$.
\end{proof}
An analogous statement holds true for ${}_{\tilde\chi}\hspace{-1pt}\wp$ with replacement
of $\A_\chi$ by ${}_{\tilde\chi}\A$, which becomes a left $\B^\Lc_\chi$-module.
\begin{propn}
\label{B-A-inv}
The projections $\wp_\chi$ and ${}_{\tilde \chi}\wp$ coincide on $\B^\Ha$.
\end{propn}
\begin{proof}
An invariant element  $a=\sum_i a_i\la_i \in \B^\Ha$ can be presented
in the form
$$
\sum_k (\Theta_1\tr\la_k) \bar\gm(\Theta_2)\tr a_k =\sum_k \bigl(\Theta_1\bar\gm^2(\Theta_2)\tr\la_k\bigr)  a_k =\sum_k (\vt\tr\la_k) a_k.
$$
This immediately implies
$\wp_\chi(a)=\sum_i a_i\chi(\la_i) ={}_{\tilde \chi}\wp(a)$.
\end{proof}
The mapping $\wp_\chi$ induces a surjective algebra homomorphism $\B^\Ha\to I_\chi$ of $\Ha$-invariants.
\section{Dynamical Galois extension}
In the present section we generalize the notion of Hopf-Galois extension to the case of dynamical algebras.
First we recall the standard definition of the Galois map for $\Ha^*$-comodule algebras.

Let $\B$ be a right $\Ha^*$-comodule algebra.
Recall that we use the symbolic notation $b^{[0]}\tp b^{\langle 1\rangle}=\delta(b)$
for right $\Ha^*$-coactions. We regard  $\B$ as a left  $\Ha$-module algebra under
the action $h\tp b\mapsto b^{[0]}\langle b^{\langle 1\rangle}, h\rangle$,
$h\in \Ha$ and $b\in \B$.
Equip the tensor product $\B\tp \Ha^*$ with the left and right $\B$-actions
$$
a.(b\tp \eta).c=abc^{[0]}\tp \eta c^{\langle 1\rangle},
$$
for $a,b,c\in\B$, $\eta\in \Ha^*$.
These actions commute with each other and make $\B\tp \Ha^*$ a $\B$-bimodule.

Let $\B^\Ha\subset \B$ denote the subalgebra of $\Ha$-invariants. Consider the map $\check\Ga\colon\B\tp \B\to \B\tp \Ha^*$ defined by
$
a\tp b\mapsto a b^{[0]}\tp b^{\langle 1\rangle}.
$
This map is a homomorphism of $\B$-bimodules, and
it is  factored through a homomorphism  $\Ga\colon\B\tp_{\B^\Ha} \B\to \B\tp \Ha^*$ called the
Galois map.
The algebra $\B$ is called Galois extension of $\B^\Ha$ if $\Ga$ is a bijection
and weak  Galois extension if $\Ga$ is a surjection.

\begin{lemma}
The vector subspace $\J_\chi=\Span\bigl(b\la\tp \eta-b\tp (\la\rl{\chi}\eta)\bigr)$, $b\in \B$, $\eta\in \Ha^*$, $\la\in \Lc$,
is a $\B$-sub-bimodule in $\B\tp \Ha^*$.
\end{lemma}
\begin{proof}
The subspace $\J_\chi$ is obviously preserved by the left $\B$-action:
$$
a.(b\la\tp \eta-b\tp \la\rl{\chi}\eta)=(ab)\la\tp \eta-(ab)\tp \la\rl{\chi}\eta,
$$
for $a,b\in \B$, $\eta\in \Ha^*$, and $\la\in \Lc$.
The following calculation demonstrates the right $\B$-invariance of $\J_\chi$:
$$
(b\la\tp \eta-b\tp \la\rl{\chi}\eta).a=b\la a^{[0]}\tp \eta a^{\langle 1\rangle}-ba^{[0]}\tp (\la \rl{\chi}\eta) a^{\langle 1\rangle}.
$$
The first summand in the right-hand side can be rewritten as
$
b(\la^{(1)}\tr a^{[0]})\la^{[\infty]}\tp \eta a^{\langle 1\rangle}=ba^{[0]}\la^{[\infty]}\tp \eta  (a^{\langle 1\rangle}\rcr\la^{(1)}),
$
while the second as
$
-ba^{[0]}\tp \la^{[\infty]} \rl{\chi}\bigl( \eta (a^{\langle 1\rangle}\rcr\la^{(1)})\bigr).
$
This proves the statement.
\end{proof}
Remark that the quotient $\B\tp \Ha^*/\J_\chi$ is isomorphic   to
the tensor product $\B \tp_\Lc \Ha^*$ as a left $\B$-module. As a vector space, it coincides with $\A\tp \Ha^*$, and we denote it by $\A\rt{\chi}\Ha^*$.
Thus the latter becomes a $\B$-bimodule. We denote by $\mathrm{P}_\chi$ the projection
$\B \tp \Ha^*\to\A\rt{\chi}\Ha^*$ along $\J_\chi$.

\begin{thm}
The map $\mathrm{P}_\chi\circ\check\Ga\colon\B\tp\B\to \A\tp_\chi \Ha^*$ factors through a $\B$-bimodule homomorphism $\Ga_\chi\colon\A_\chi\tp_{I_\chi} {}_{\tilde\chi}\A\to \A\rt{\chi}\Ha^*$.
\end{thm}
\begin{proof}
First of all, the map $\mathrm{P}_\chi\circ\check\Ga$ induces a $\B$-bimodule map $\check\Ga_\chi\colon\A_\chi\tp {}_{\tilde\chi}\A\to \A\rt{\chi}\Ha^*$.
This follows from the equality
$(\mathrm{P}_\chi\circ\check\Ga)(\la\tp \mu)=\chi(\la)\tilde\chi(\mu)(1\tp 1)$  which holds true  for all $\la\in \Lc$ by the very construction
and for all $\mu\in \Lc$ by Lemma \ref{scalar}.
Further, for all $a,b\in \A$ and $c\in I_\chi$ one has
$$
\check\Ga_\chi(a\tp c\rd{\tilde\chi} b)=\check\Ga_\chi(a\tp c).b=(a\st{1}c^{[0]}\tp \rho_1\rl{\chi}c^{\langle 1\rangle}).b=\check\Ga_\chi(a\st{\chi}c\tp 1).b=\check\Ga_\chi(a\st{\chi}c\tp b).
$$
This  implies that the map $\check\Ga_\chi$ factors through a map $\Ga_\chi\colon \A_\chi\tp_{I_\chi} {}_{\tilde\chi}\A\to \A\rt{\chi}\Ha^*$,
which is a homomorphism of $\B$-bimodules.
\end{proof}
\begin{definition}
\label{def_Galois}
We call $\Ga_\chi$  the  $\chi$-Galois map.
The dynamical algebra
$\A$ is called a (weak) $\chi$-Galois extension of $I_\chi$ if
$\Ga_\chi$ is an (epimorphism) isomorphism.
\end{definition}
Definition \ref{def_Galois} reduces to the standard Galois map
for an ordinary $\Ha$-module algebra $\A$ upon setting $\Lc=\Ha$ and $\chi=\ve$.
Applying our construction to this situation we conclude that the standard Galois map (which is independent on $\chi$) is a morphism
of not only  $A$-bimodules, but of $A\rtimes\Ha$-bimodules too. The two-sided
action of $A\rtimes\Ha$ on $\A\tp_{\A^\Ha} \A$ is straightforward while
it is not obvious on  $\A\tp \Ha^*$. For instance, the left $\Ha$-module $\A\tp \Ha^*$
is a tensor product of the default module $\A$ and the module $\Ha^*$
under the action $h.\al= \al \rcr \gm(h)$. The right $\Ha$-module
$\A\tp \Ha^*$ is a tensor product of the trivial module $\A$ and the module $\Ha^*$
with the right action $\al.h= \sum_{i,j}\bigl(\al \eta^i\gm(\eta^j)\bigr)\rcr \gm(h_ih h_j)$.
Here $\{h_i\}\subset \Ha$ is a basis and $\{\eta^i\}\subset \Ha^*$ is its dual.
Note that for either commutative or co-commutative $\Ha$ this
action is
expressed through the left co-regular action  and antipode  by $\al.h= \gm(h)\lcr\al$.

Below we  display the construction of Galois maps on a commutative diagram
of $\B$-bimodules.
\be
\begin{array}{c}
\begin{picture}(230,125)
\put(5,110){$\B\tp \B$}
\put(200,110){$\B\tp \Ha^*$}
\put(40,113){\vector(1,0){150}}
\put(22,105){\vector(0,-1){50}}
\put(112,72){\vector(0,-1){50}}

\put(25,106){\vector(3,-1){68}}
\put(147,84){\vector(3,1){65}}
\put(95,78){$\B\tp_{\B^\Ha}\B$}
\put(217,105){\vector(0,-1){50}}

\put(0,43){$\A_\chi\tp {}_{\tilde\chi}\A$}
\put(200,43){$\A\tp_\chi \Ha^*$}
\put(50,48){\vector(1,0){145}}

\put(24,40){\vector(3,-1){65}}
\put(150,18){\vector(3,1){60}}
\put(91,13){$\A_\chi\tp_{I_\chi} {}_{\tilde\chi}\A$}

\put(-20,75){$\wp_\chi\tp {}_{\tilde\chi}\hspace{-1.5pt}\wp$}

\put(125,117){$\check\Ga$}
\put(125,52){$\check\Ga_\chi$}
\put(220,75){$\mathrm{P}_\chi$}
\put(175,80){$\Ga$}
\put(175,12){$\Ga_\chi$}
\end{picture}
\end{array}
\ee
The front downward arrow is due to Proposition \ref{B-A-inv}.
Note that $\A/I_\chi$ is a weak $\chi$-Galois extension if so is  $\B/\B^\Ha$, as follows from the diagram.

Recall from Lemma \ref{l-r} that the subspace of $\chi$-generators in $\A_\chi$ is equal to
the subspace of  $\tilde \chi$-generators in ${}_{\tilde\chi}\A$. By Theorem \ref{End},
the operations $\st{\chi}$ and $\rd{\tilde\chi}$ coincide when restricted to that
subspace and define the algebra $\A^\Lc_\chi$.
\begin{propn}
Restriction of the $\Ga_\chi$ to $\A^\Lc_\chi\tp_{I_\chi} \tp \A^\Lc_{\chi}$
is the ordinary Galois map and takes values in $\A^\Lc_\chi\tp (\Ha^\chi)^*\subset \A^\Lc_\chi\tp \Ha^*$.
\end{propn}
\begin{proof}
The first statement follows from the equalities
$$
a\st{1}b^{[0]}\tp \rho_{1}\btr b^{\langle 1\rangle}=a\st{1}b^{[0]}\tp \chi(\rho_{1}) b^{\langle 1\rangle}=a\st{\chi}b^{[0]}\tp  b^{\langle 1\rangle}
$$
valid for all $a,b\in \A^\Lc_\chi$. The coaction $\A^\Lc_\chi\to \A^\Lc_\chi\tp \Ha^*$
takes values in $\A^\Lc_\chi\tp (\Ha^\chi)^*$, as the $\Ha$-action on $\A^\Lc_\chi$ factors through an action of $\Ha^\chi$.
\end{proof}

The Galois map is independent on the choice of base algebra in the following sense.
Suppose we are given  a homomorphism of bases $(\Ha,\Lc)\stackrel{\id,\varpi}{\longrightarrow} (\Ha,\Lc')$ identical on $\Ha$.
Given a dynamical algebra $\A$ over the base $(\Ha,\Lc)$ we make it a dynamical algebra over $(\Ha,\Lc')$, with the multiplication
$$
\A\tp \A\to \A\tp \Lc\stackrel{\id \tp \varpi}{\longrightarrow} A\tp \Lc'.
$$
This amounts to an $\Ha$-algebra map $\A\rtimes \Lc\to \A\rtimes \Lc'$ identical on $\A$ and given by $\varpi$ on $\Lc$. A
character $\chi'\in \hat \Lc'$ gives rise to the character $\chi:=\chi'\circ \varpi\in \hat \Lc$.
We have also $\tilde \chi:=\tilde \chi'\circ \varpi$, as  $\varpi$ commutes with the
action of $\vt\in \D(\Ha)$ on $\Lc$ and $\Lc'$ (because the base homomorphism is identical on $\Ha$).
It is obvious that $\Ga_{\chi}$ constructed via $\A\rtimes \Lc$ is equal to $\Ga_{\chi'}$ constructed via $\A\rtimes \Lc'$.
Below we illustrate this on the example $\Lc'=\Ha$.

We fix an invariant character $\chi\in \hat\Lc$ and consider the homomorphism $\iota_\chi\colon\Lc\to \Ha$.
As was mentioned earlier, this amounts to
a  homomorphism of dynamical bases $(\Ha,\Lc)\to (\Ha,\Ha)$ identical on $\Ha$.
The group $\hat\Ha \subset \Ha^*$ acts on the algebra  $\Ha$ by base algebra automorphisms through the left regular action
$\varpi_\al(h):=\al\lcr  h=h^{(1)} \langle\al,h^{(2)}\rangle$, see Example \ref{al-auto}. Under this action,
 $\al\circ \iota_\chi= \iota_{\al\cdot\chi}$.
Denote by $A\rtimes_\chi\! \Ha$ the dynamical smash product of $\A$ and $\Ha$.
For  $\al\in \hat \Ha$ the algebras $A\rtimes_\chi \! \Ha$
and $A\rtimes_{\al\cdot\chi} \!\Ha$ are isomorphic, with the isomorphism
identical on $\A$ and given by $\varpi_\al$ on $\Ha$.

All the characters of the base algebra $\Ha$ are invariant and constitute the group $\hat\Ha$.
For all $\al\in \hat \Ha$ we construct the left and right $A\rtimes_\chi\! \Ha$-induced modules $\A_\al$ and ${}_{\al}\A$ as in
Section \ref{secLBI}.
The homomorphism $\B\to A\rtimes_\chi\! \Ha$ makes  $\A_\al$ and ${}_{\al}\A$ left and right $\B$-modules, which
are isomorphic to $\A_{\al\cdot \chi}$ and ${}_{\al\cdot \chi}\A$, respectively.
The subspace $I_\al\in \A$ is formed by $\Ha$-invariants.
We have the following isomorphisms
$(I_{\al\cdot\chi})_{op}\simeq (I_\al)_{op} \simeq\End_{A\rtimes_\chi\hspace{-0.5pt}  \Ha}(\A_{\al})$,
where the right isomorphism is due to Theorem \ref{End} applied to base algebra $\Lc=\Ha$.

The homomorphism $\B\to \B':=A\rtimes_\chi \Ha$  gives rise to the commutative diagrams
$$
\begin{diagram}
   \dgARROWLENGTH=0.5\dgARROWLENGTH
\node{\B \tp \B}\arrow{e,t}{\check\Ga}\arrow{s,t}{}\node{\B \tp \Ha^*}\arrow{s,t}{}
\node{\A_{\al\cdot\chi} \tp_{I_\chi}{}_{\al\cdot\tilde\chi}\A}\arrow{e,t}{\Ga_{\al\cdot\chi}}\arrow{s,r}{\wr}\node{\A \tp \Ha^*}\arrow{s,r}{\wr}
\\
\node{\B' \tp \B'}\arrow{e,t}{\check\Ga'}\node{\B' \tp \Ha^*}
\node{\A_\al \tp_{I_\al} {}_{\tilde\al}\A}\arrow{e,t}{\Ga'_\al}\node{\A \tp \Ha^*}
\end{diagram}
,\quad \forall \al\in \hat \Ha.
$$
In the right square we have used the equality $\al\cdot\tilde\chi=\widetilde{\al\cdot\chi}$, cf. Section \ref{SecICh}.
These diagrams illustrate independence of the Galois map on the choice of base algebra.

\section{Morita context for dynamical smash product}
In the present section we construct  a Morita context for the dynamical smash product
with any invariant character $\chi\in \hat\Lc$ and a left generator $e\in \Lc$.
Recall  that two algebras are called Morita equivalent if their module categories are equivalent.
Morita context is a weaker relation than equivalence.
Given two algebras $C$ and $\B$, a Morita context is a quadruple
$
\bigl(P,Q,[.,.],(.,.)\bigr)
$
consisting of an $C-\B$ bimodule $P$, an $\B-C$ bimodule $Q$ and two maps
$$
[.,.]\colon P\tp_\B Q\to C,
\quad
(.,.)\colon Q\tp_C P\to \B
$$
obeying the following conditions.
Firstly, they are, respectively, $C$- and $\B$-bimodule maps.
Secondly, for all $p,p'\in P$ and $q,q'\in Q$ they satisfy the "associativity"  condition
$$
p'\dashv(q,p)=[p',q]\vdash p,
\quad
q'\dashv[p,q]=(q',p)\vdash q.
$$
Here $\vdash$ and $\dashv$ denote the left and right actions of $C$ and $\B$.
Morita context defines two functors between categories of $C$-and $\B$-modules and
becomes equivalence if $[.,.]$ and $(.,.)$ are surjective.

We construct a Morita context for the algebras $\B$ and $\A_\chi^\Lc$ under the assumption
that the subspace $\Lc^\chi$ of $\chi$-generators is not zero. The construction depends
on a generator $e\in \Lc^\chi$.
\begin{lemma}
For every $e\in \Lc$ the map $\check\Pi_e\colon\B\tp \Ha^*\to \B$, $a\tp \eta\mapsto a \langle \eta, e^{(1)}\rangle e^{[\infty]}$
is a homomorphism of $\B$-bimodules. If $e$ is a $\chi$-generator, then this
map factors through a map $\Pi_e\colon\A\tp_\chi \Ha^*\to \B$.
\end{lemma}
\begin{proof}
The statement is obvious as to the left $\B$-action. Let us demonstrate that $\check\Pi_e$  is a right $\B$-homomorphism.
For all $a,b\in \B$ and $\eta\in \Ha^*$ the element $ab^{[0]}\tp \eta b^{\langle 1 \rangle}$ goes to
$$
a b^{[0]} \langle \eta b^{\langle 1 \rangle}, e^{(1)}\rangle e^{[\infty]}
=a b^{[0]} \langle \eta, e^{(1)}\rangle\langle b^{\langle 1 \rangle}, e^{(2)}\rangle e^{[\infty]}
=a (e^{(2)}\tr  b)\langle \eta, e^{(1)}\rangle e^{[\infty]}
=a \langle \eta, e^{(1)}\rangle e^{[\infty]} b,
$$
as required.

Further, suppose that $e$ is a left $\chi$-generator. Then
$
 \langle \eta, e^{(1)}\rangle\la   e^{[\infty]}   =
\langle \la\rl{\chi} \eta, e^{(1)}\rangle e^{[\infty]}
$
    for all $\la\in \Lc$ and $\eta\in \Ha^*$. Hence
$$
\check\Pi_e\colon a\la\tp \eta-a\tp (\la\rl{\chi} \eta)\mapsto a \la   \langle \eta, e^{(1)}\rangle e^{[\infty]}-a \langle \la\rl{\chi} \eta, e^{(1)}\rangle e^{[\infty]}=0.
$$
Thus  $\check\Pi_e(\J_\chi)=\{0\}$, and $\check\Pi_e$ factors through a $\B$-bimodule map $\Pi_e\colon\A\tp_\chi \Ha^*\to \B$.
This completes the proof.
\end{proof}

We call a left $\chi$-generator $e$ cyclic if $\Lc$ is a cyclic right $\Ha^*$-module generated by $e$.
In particular, in the special case $\Lc=\Ha$ every character of $\Ha$ is cyclic. This follows
from the fact that the map $\eta\mapsto t\rcr \eta$, where $t\not= 0$ is a left integral
of $\Ha$, is a linear isomorphism.
\begin{lemma}
Suppose that $\A$ is a weak $\chi$-Galois extension. If the left $\chi$-generator $e$ is cyclic,
then the map $\A_\chi\tp_{I_\chi} {}_{\tilde\chi}\A\to \A\tp_\chi \Ha^*\stackrel{\Pi_{e_\chi}}{\longrightarrow} \B$
induces an epimorphism $\A_\chi\tp_{\A^\Lc_\chi} {}_{\tilde\chi}\A\to \B$
given by the assignment $a\tp b\mapsto a e b$.
\end{lemma}
\begin{proof}
The map in question is given by the assignment
$a\tp b\mapsto a (e^{(1)}\tr b)e^{[\infty]}=ae b$. This implies the statement.
\end{proof}

Put
$
P:=\A_\chi, $ and $Q:={}_{\tilde\chi}\A.$
By construction $P$ and $Q$ are, respectively, $\B-\A^\Lc_\chi$, and $\A^\Lc_\chi-\B$ bimodules.
The $\B$-module structures are obtained by induction from the $\Ha$-module $\Bbbk_\chi$.
Define the maps
$$
[.,.]\colon P\tp_{\A^\Lc_\chi} Q\to \B
,
\quad
(.,.)\colon Q\tp_{\B} P\to \A_\chi^\Lc
$$
by setting
$$
[a ,b]=a e b
,
\quad
(a,b)=\iota_\chi(e)\tr(a\st{\chi}b).
$$
\begin{propn}
The quadruple $\bigl(P,Q,[.,.],(.,.)\bigr)$ is a Morita context
for $\B$ and $\A^\Lc_\chi$.
It becomes equivalence if the extension $\A/I_\chi$ is weak Galois,  the generator $e$ is cyclic,
and $\A^\Lc_\chi=\iota_\chi(e)(\A)$.
The similar statement is true if  $\B$ is replaced by $\A\rtimes_\chi \Ha$ and $\A^\Lc_\chi$ by $I_\chi$.
\end{propn}
\begin{proof}
It is straightforward to check that $\bigl(P,Q,[.,.],(.,.)\bigr)$ defined as above is
a Morita context indeed.
The map $(.,.)$ is an epimorphism if $\A^\Lc=\iota_\chi(e)(\A)$. Therefore the Morita context becomes equivalence if and only
if the map $[.,.]$ is surjective. That is the case if $\A$ is a weak $\chi$-Galois extension and $e$ is cyclic.
Setting $\Lc=\Ha$ we get $\A^\Lc_\chi=I_\chi$, so the last statement follows from the previous two.
\end{proof}

The constructions of this section take an especially simple form if the generator $e=e_\chi$ is
an idempotent; recall that $\tilde \chi=\chi$ in this case.
For the algebra $\B=\A\rtimes \Lc$,  its left induced module $\A_\chi$ is realized as the left ideal
$\B e_\chi$, while the right induced module ${}_{\chi}\A$ as the
right ideal $e_\chi \B$.
The algebra $\A^\Lc_\chi$ is isomorphic to the
subring $e_\chi \B e_\chi \subset \B$.
The Galois map can be presented as $a\tp b\mapsto a e_\chi b^{[0]}\tp b^{\langle 1\rangle}$.
The Morita context
is defined through the maps
$$
[a e_\chi,e_\chi b]=a e_\chi b
,
\quad
(e_\chi a,b e_\chi)=e_\chi(ab)e_\chi.
$$
where $P=\B e_\chi=\A e_\chi $ and $Q=e_\chi\B=e_\chi \A.$
\section{Ring theoretical characterization of $\chi$-Galois extension}
In this last section we generalize Theorem 1.2 of \cite{CFM} characterizing the standard Hopf-Galois extension,
for the case of dynamical algebras.
\begin{lemma}[Frobenius reciprocity]
Suppose $\B$ is a dynamical smash product $\A\rtimes \Lc$.
Let $M$ be a left $\B$-module and fix a character $\chi\in \hat \Lc$.
Then $\Hom_\Lc(\Bbbk_\chi,M)\simeq \Hom_\B(\A_\chi,M)$ as vector spaces.
The analogous statement is true for right modules.
\end{lemma}
This lemma is analogous to Lemma 0.6 of \cite{CFM} and, upon replacement of $\Lc$ by $\Ha$,  plays the same role in the following theorem
for  $\B=\A\rtimes_\chi\Ha$.
\begin{thm}
Let $\Lc$ be a base algebra over a finite dimensional Hopf algebra $\Ha$ and let $\A$ be
a dynamical $(\Ha,\Lc)$-algebra. Fix an invariant character $\chi\in \hat \Lc$. Then the following are equivalent
\begin{enumerate}
\item[(1)] $\A_\chi/I_\chi$ is right weak $\Ha^*$-Galois.
\item[(2)]
\begin{enumerate}
\item[(a)] The map $\A\rtimes_\chi\Ha\to \End
_{I_\chi}(\A_\chi)$ is an algebra isomorphism,
\item[(b)] $\A_\chi$ is a finitely generated projective right $I_\chi$-module.
\end{enumerate}
\item[(3)] $\A_\chi$ is a left $\A\rtimes_\chi\Ha$-generator.
\item[(4)] if $0\not=t\in \int_l\subset\Ha$, then the map $[.,.]\colon\A_\chi\tp_{I_\chi} {}_{\tilde\chi}\A\to \A\rtimes_\chi\Ha$
given by $[a,b]= a t b$ is surjective.
\item[(5)] For any left $\A\rtimes_\chi\Ha$-module the map $\A_\chi\tp_{I_\chi} M^\Ha\to M$
given by $a\tp m\mapsto a m$ is a left $\A\rtimes_\chi\Ha$-module isomorphism.
\end{enumerate}
\end{thm}
The proof of this theorem is literally the same as the proof of Theorem 1.2 in \cite{CFM}.
The only difference is the argumentation for the equality $(\A\rtimes_\chi\Ha)^\Ha=t\A$ (here
$\A\rtimes_\chi\Ha$ is understood as a left regular $\Ha$-module)
in the proof of the implication $(5)\Rightarrow (4)$.
The inclusion $(\A\rtimes_\chi\Ha)^\Ha\subset t\A$  follows from $\dim \int_l=1$ and from the fact that the left $\Ha$-module $\A\rtimes_\chi\Ha$ is freely generated
by $\A$; the assumption for $\A$ being an algebra is unnecessary.

\end{document}